\documentclass[12pt]{amsart}
\usepackage{amssymb} 
\usepackage{a4, amsthm} 

\theoremstyle{plain}
\newtheorem{theorem}{Theorem}[section]

\newtheorem{prop}[theorem]{Proposition}

\theoremstyle{definition}

\font\ninerm=cmr9  \font\eightrm=cmr8  \font\sixrm=cmr6
\font\ninei=cmmi9  \font\eighti=cmmi8  \font\sixi=cmmi6
\font\ninesy=cmsy9 \font\eightsy=cmsy8 \font\sixsy=cmsy6
\font\ninebf=cmbx9 \font\eightbf=cmbx8 \font\sixbf=cmbx6
\font\nineit=cmti9 \font\eightit=cmti8 
\font\ninett=cmtt9 \font\eighttt=cmtt8 
\font\ninesl=cmsl9 \font\eightsl=cmsl8

\font\twelverm=cmr12 at 15pt
\font\twelvei=cmmi12 at 15pt
\font\twelvesy=cmsy10 at 15pt
\font\twelvebf=cmbx12 at 15pt
\font\twelveit=cmti12 at 15pt
\font\twelvett=cmtt12 at 15pt
\font\twelvesl=cmsl12 at 15pt
\font\twelvegoth=eufm10 at 15pt

\font\tengoth=eufm10  \font\ninegoth=eufm9
\font\eightgoth=eufm8 \font\sevengoth=eufm7 
\font\sixgoth=eufm6   \font\fivegoth=eufm5
\newfam\gothfam \def\goth{\fam\gothfam\tengoth} 
\textfont\gothfam=\tengoth
\scriptfont\gothfam=\sevengoth 
\scriptscriptfont\gothfam=\fivegoth

\catcode`@=11
\newskip\ttglue

\def\tenpoint{\def\rm{\fam0\tenrm}
  \textfont0=\tenrm \scriptfont0=\sevenrm
  \scriptscriptfont0\fiverm
  \textfont1=\teni \scriptfont1=\seveni
  \scriptscriptfont1\fivei 
  \textfont2=\tensy \scriptfont2=\sevensy
  \scriptscriptfont2\fivesy 
  \textfont3=\tenex \scriptfont3=\tenex
  \scriptscriptfont3\tenex 
  \textfont\itfam=\tenit\def\it{\fam\itfam\tenit}%
  \textfont\slfam=\tensl\def\sl{\fam\slfam\tensl}%
  \textfont\ttfam=\tentt\def\tt{\fam\ttfam\tentt}%
  \textfont\gothfam=\tengoth\scriptfont\gothfam=\sevengoth 
  \scriptscriptfont\gothfam=\fivegoth
  \def\goth{\fam\gothfam\tengoth}
  \textfont\bffam=\tenbf\scriptfont\bffam=\sevenbf
  \scriptscriptfont\bffam=\fivebf
  \def\bf{\fam\bffam\tenbf}%
  \tt\ttglue=.5em plus.25em minus.15em
  \normalbaselineskip=12pt \setbox\strutbox\hbox{\vrule
  height8.5pt depth3.5pt width0pt}%
  \let\big=\tenbig\normalbaselines\rm}

\def\ninepoint{\def\rm{\fam0\ninerm}
  \textfont0=\ninerm \scriptfont0=\sixrm
  \scriptscriptfont0\fiverm
  \textfont1=\ninei \scriptfont1=\sixi
  \scriptscriptfont1\fivei 
  \textfont2=\ninesy \scriptfont2=\sixsy
  \scriptscriptfont2\fivesy 
  \textfont3=\tenex \scriptfont3=\tenex
  \scriptscriptfont3\tenex 
  \textfont\itfam=\nineit\def\it{\fam\itfam\nineit}%
  \textfont\slfam=\ninesl\def\sl{\fam\slfam\ninesl}%
  \textfont\ttfam=\ninett\def\tt{\fam\ttfam\ninett}%
  \textfont\gothfam=\ninegoth\scriptfont\gothfam=\sixgoth 
  \scriptscriptfont\gothfam=\fivegoth
  \def\goth{\fam\gothfam\tengoth}
  \textfont\bffam=\ninebf\scriptfont\bffam=\sixbf
  \scriptscriptfont\bffam=\fivebf
  \def\bf{\fam\bffam\ninebf}%
  \tt\ttglue=.5em plus.25em minus.15em
  \normalbaselineskip=11pt \setbox\strutbox\hbox{\vrule
  height8pt depth3pt width0pt}%
  \let\big=\ninebig\normalbaselines\rm}

\def\eightpoint{\def\rm{\fam0\eightrm}
  \textfont0=\eightrm \scriptfont0=\sixrm
  \scriptscriptfont0\fiverm
  \textfont1=\eighti \scriptfont1=\sixi
  \scriptscriptfont1\fivei 
  \textfont2=\eightsy \scriptfont2=\sixsy
  \scriptscriptfont2\fivesy 
  \textfont3=\tenex \scriptfont3=\tenex
  \scriptscriptfont3\tenex 
  \textfont\itfam=\eightit\def\it{\fam\itfam\eightit}%
  \textfont\slfam=\eightsl\def\sl{\fam\slfam\eightsl}%
  \textfont\ttfam=\eighttt\def\tt{\fam\ttfam\eighttt}%
  \textfont\gothfam=\eightgoth\scriptfont\gothfam=\sixgoth 
  \scriptscriptfont\gothfam=\fivegoth
  \def\goth{\fam\gothfam\tengoth}
  \textfont\bffam=\eightbf\scriptfont\bffam=\sixbf
  \scriptscriptfont\bffam=\fivebf
  \def\bf{\fam\bffam\eightbf}%
  \tt\ttglue=.5em plus.25em minus.15em
  \normalbaselineskip=9pt \setbox\strutbox\hbox{\vrule
  height7pt depth2pt width0pt}%
  \let\big=\eightbig\normalbaselines\rm}

\def\twelvepoint{\def\rm{\fam0\twelverm}
  \textfont0=\twelverm\scriptfont0=\tenrm
  \scriptscriptfont0\sevenrm
  \textfont1=\twelvei\scriptfont1=\teni
  \scriptscriptfont1\seveni 
  \textfont2=\twelvesy\scriptfont2=\tensy
  \scriptscriptfont2\sevensy 
   \textfont\itfam=\twelveit\def\it{\fam\itfam\twelveit}%
  \textfont\slfam=\twelvesl\def\sl{\fam\slfam\twelvesl}%
  \textfont\ttfam=\twelvett\def\tt{\fam\ttfam\twelvett}%
  \textfont\gothfam=\twelvegoth\scriptfont\gothfam=\ninegoth 
  \scriptscriptfont\gothfam=\sevengoth
  \def\goth{\fam\gothfam\twelvegoth}
  \textfont\bffam=\twelvebf\scriptfont\bffam=\ninebf
  \scriptscriptfont\bffam=\sevenbf
  \def\bf{\fam\bffam\twelvebf}%
  \tt\ttglue=.5em plus.25em minus.15em
  \normalbaselineskip=12pt \setbox\strutbox\hbox{\vrule
  height9pt depth4pt width0pt}%
  \let\big=\twelvebig\normalbaselines\rm}

\def\C{{\mathbb{C}}}

\def\P{{\mathbb{P}}}

\def\nsmallskip{\smallskip\noindent}
\def\bbigskip{\bigskip\bigskip}
\def\nbigskip{\bigskip\noindent}

\def\nmedskip{\medskip\noindent}
\def\buildunder#1#2{\mathrel{\mathop{\kern0pt #2}
\limits_{#1}}}

\def\qq{/\kern-.185em /}

\def\mn{\medskip\noindent}
\def\bn{\bigskip\noindent}
\def\REM #1{}

\begin{document}

\title[Hyperbolic manifolds]{Hyperbolic manifolds whose envelopes of holomorphy are
not hyperbolic }

\bbigskip

\author[Geatti]{L. Geatti}
\author[Iannuzzi]{A. Iannuzzi}
\author[Loeb]{J.-J. Loeb}
\address{Laura Geatti and Andrea Iannuzzi: Dip.~di Matematica,
II Universit\`a di Roma  ``Tor Vergata", Via della Ricerca Scientifica,
I-00133 Roma, Italy.} 
\email{geatti@mat.uniroma2.it,
       iannuzzi@mat.uniroma2.it}
\address{Jean-Jacques Loeb, D\'epartement de Math\'ematiques et Informatiques, Universit\'e d'Angers, 2 Boulevard Lavoisier, 49045 Angers cedex 01}
\email{Jean-Jacques.Loeb@univ-angers.fr}
  
\thanks {\ \ {\it Mathematics Subject Classification (2000):} 32Q45, 32D10, 32M05}

\thanks {\ \ {\it Key words}:
hyperbolic manifold,  envelope of holomorphy, symmetric space}

\begin{abstract}
We present a family of examples of two dimensional, hyperbolic complex
manifolds whose envelopes of holomorphy are not hyperbolic.  
\end{abstract}
 
\maketitle

\bigskip
\bbigskip

\medskip
In this note we present 
a family of hyperbolic complex manifolds whose envelopes of holomorphy are
 not hyperbolic. Here we consider hyperbolicity  in the sense of Kobayashi.
In Isaev's classification \cite {Is} 
 of two-dimensional hyperbolic manifolds
 with three-dimensional automorphism group,  these  manifolds appear  
 as subdomains of $\,\Delta \times \P^1$, with $\,\Delta\,$  the unit disk in $\,\C\,$
 and $\,\P^1\,$  the one-dimensional complex projective  space.
 Here we consider their
 realizations as $\,SU(1,1)$-invariant domains in the Lie group
 complexification $\,SL(2,\C)/U(1)^\C\,$ of  the symmetric space $\,SU(1,1) / U(1)$. Then, by applying
 the  univalence results for Stein, equivariant Riemann domains
 over Lie group complexifications of rank-one, Riemannian 
 symmetric spaces  obtained in \cite{GeIa},  
 we can explicitely determine their envelopes of holomorphy.   Such envelopes   turn out to be all biholomorphic to $\,\Delta \times \C$. In particular,  
 they are not hyperbolic.

 
 \mn
 
Let  $\,G\,$ be the Lie group $\,SU(1,1)\,$. Consider the holomorphic action  of its universal complexification $\,G^\C=SL(2,\C) \,$  
on $\,\P^1\times \P^1\,$  defined  by
\begin{equation}
\label{ACTION}
g\cdot([z_1:z_2],[w_1:w_2]):=(\,g\cdot [z_1:z_2],\  \overline{\sigma(g)}\cdot  [w_1:w_2] \,)\,,
\end{equation}

\nsmallskip
where $\,\sigma(g)=I_{1,1}\,^t \bar g^{-1}I_{1,1}\,$ denotes the conjugation of $\,G^\C\,$ relative to $\,G\,$  
(here $\,I_{1,1}\,$ the diagonal matrix representing the standard 
hermitian form of signature $\,(1,1)$).
 Denote by $\,K\,$ a maximal compact subgroup of $\,G\,$ and by $\,K^\C\,$ its complexification.
  The quotient $\,G^\C/K^\C\,$ can be  identified with the unique open 
  $\,G^\C$-orbit in  $\,\P^1 \times \P^1\,$ given by
$$  G^\C \cdot ([1:0],[1:0])
=\{\, ([z_1:z_2],[w_1:w_2]) \in \P^1 \times \P^1  \ : \  z_1w_1- z_2 w_2\not=0\,\}\,.$$

\nbigskip
A global slice for the  $G$-action  on $\,G^\C/K^\C\,$ by left translations is represented by the following diagram (cf. \cite{GeIa}, Sect.~4.1) 
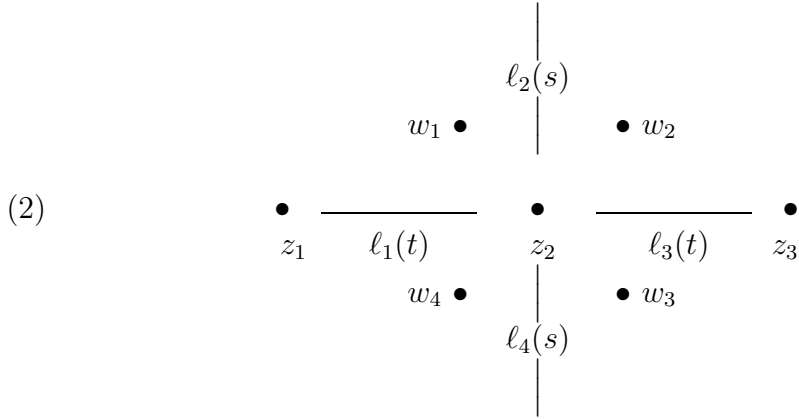
\begin{equation}
\label{DIAGRAM1}
\begin{matrix}  
&  &\Big|& &   \cr 
&\quad \quad  \quad \quad \quad \quad \quad \quad &   \ell_2( s) & &  \cr
&\ \ \ \  \quad \quad \quad \quad w_1 \ \bullet & \Big|&\bullet \  w_2 \quad \quad \quad \quad  \quad \ \ &  \cr 
&  && &   \cr 
 &\ \ \ \   \bullet \quad \overline{ \quad \quad \quad \quad \quad  }&\bullet&\overline{ \quad \quad \quad \quad \quad   }\quad \bullet  \ \ \ \ \ \ \ &\cr  
  &\  z_1 \quad \ \ \ \ell_1(t)  \      &\  z_2 & \ell_3(t)\ \ \   \quad  z_3\ \ \cr
  & \ \ \ \  \quad \quad \quad  \quad w_4\ \bullet&\Big|&\bullet\ w_3\ \  \quad \quad \quad \quad  \quad&  \cr 
  & &  \ell_4(s)& &  \cr 
  &  &\Big|& &   \cr 
  \end{matrix}
  \end{equation}

\bn
All elements in the diagram, except for $\,z_1$, $\,z_2\,$ and $\,z_3$,
lie on
hypersurface
$\,G$-orbits. 
 The points 
 $\, z_1:=([1:0], [ 1:0])\,$ and $\, z_3:=([ 0:1], [ 0:1])\,$ lie on $\,G$-orbits  diffeomorphic to the symmetric space $\,G/K$.  
The point   $z_2:=([ 1 :i], [ 1:i ])\,$ lies on a  $G$-orbit  diffeomorphic to a pseudo-Riemmanian symmetric space of the same dimension as $G/K$.
The  slices $\,\ell_1, \ldots,\ell_4\,$ are defined by 
$$\ell_1(t)=([\cos {\pi\over 4}(1-t): i\sin {\pi\over 4}(1-t)], [ \cos {\pi\over 4}(1-t): i\sin {\pi\over 4}(1-t)]) \,,$$  
$$\ell_3(t)= ([ \cos {\pi\over 4}(1+t):i\sin {\pi\over 4}(1+t)], [ \cos {\pi\over 4}(1+t): i\sin {\pi\over 4}(1+t)])\,,$$
for $0<t<1$, and by 
 $$ \ell_2(s)=([  e^{s} :ie^{ -s} ], [  e^{ -s} :ie^{s} ])\,, $$
 $$ \ell_4(s)=([ e^{-s} :ie^{ s} ], [  e^{ s} :ie^{ -s} ])\,, $$ 
for $\,s>0$.
 Note that $z_2$ is the limit point of $\ell_1$, $\ell_2$, $\ell_3$  and $\ell_4$ for values of the parameters approaching zero. Similarly  one has $\ell_1(1)=z_1$ and $\ell_3(1)=z_3$.
   The points
$$\,w_1:=([ 1:0], [ 1:i]),\ \    w_2:=( [1:i],[ 0:1])\,,$$
$$\ \  w_3:=([ 1:i], [ 1:0]) \ \ w_4:=( [ 0:1], [1:i])\,.$$
represent the  four non-closed hypersurface  $\,G$-orbits containing  $\,G\cdot z_2\,$ in their closure.

\bn 
Consider the family of $\,G$-invariant domains in $\,G^\C/K^\C\,$ defined by
$$D_\beta\,= \, G\cdot (\,z_1\cup\ell_1((0,1)) \cup w_1 \cup \ell_2((0,\beta))\,) \qquad {\rm for} \quad 0<\beta<\infty \,  .$$
By the classification of Stein, $\,G$-invariant domains in $\,G^\C/K^\C\,$ (cf.~Thm.~6.1
in \cite{GeIa}), none of the   domains $\,D_\beta\,$ is Stein.  Each of them contains
a unique Levi-flat orbit  given by  

\nsmallskip
 $$\,G \cdot w_1= \{\, ([1:z],[1:w]) \in \P^1 \times \P^1  \ : \ z \in \Delta, \ \  w \in
 \partial \Delta \,\}\, \cong \Delta \times \partial \Delta\, ,$$
 \nmedskip
and a  unique totally-real orbit  given by  $$\, G \cdot z_1 =\{\, ([1:z],[1: \bar z]) \in\P^1 \times \P^1  \ : \ z \in \Delta\,\}.$$

\nsmallskip
Denote by $\,W\,$ the limit domain  
$$W=G\cdot (\, z_1\cup\ell_1(0,1)  \cup w_1 \cup \ell_2(0,\infty)\,)\, .$$
By construction $\,W\,$ contains every domain~$\,D_\beta$. One can check that $\,W\,$  coincides 
with the  domain 
 $$ \{\, ([1:z],[w_1:w_2]) \in \P^1\times \P^1  \ : \ z \in \Delta,~ w_1-zw_2\not=0\}\,.$$ 
 
\nsmallskip
It follows that $\,W\,$ is biholomorphic to $\,\Delta \times \C\,$ via the map
 $$\Delta \times \C \to W, \quad (u,v) \to ([1:u],[1+  u v:v])\, .$$
In particular, it is Stein and  is not hyperbolic.
\nsmallskip

 \nmedskip
 \begin{prop}
 For every $\,0<\beta<\infty\,$  the  $\,G$-invariant domain 
  $\,D_\beta\,$ is hyperbolic. The envelope of holomorphy of $\,D_\beta\,$ is the domain $\, W $, which is  not hyperbolic.
\end{prop}

\nmedskip
 \begin{proof}
 The fact that  $\,D_\beta\,$ is hyperbolic follows from Isaev's classification of two-dimensional hyperbolic manifolds with three dimensional automorphism group. 
Indeed note that for $\,g \in G\,$ one has  $\, \overline{\sigma(g)} = \overline g$. Thus
 the restriction to $\,G\,$ of the  $\,G^\C$-action defined in (\ref{ACTION}) agrees with 
 the $\,G$-action  given  in \cite{Is}, p.$\,$22. As a consequence, 
 $\,D_\beta\,$ coincides with an element of the family denoted there by  
 $$\widehat {\mathfrak D} ^{(1)}_t  \quad {\rm for } \quad 1< t < \infty $$
and it is hyperbolic.

 Denote by $\,E(D_\beta)\,$ the  envelope of holomorphy of $\,D_\beta$. Since $\,G^\C/K^\C\,$
is Stein, the inclusion of $\,D_\beta\,$ in $\,G^\C/K^\C\,$ extends
to a local biholomorphism $\,p\, $ from $\, E(D_\beta) \,$ to  $\, G^\C/K^\C \,$ (cf.~\cite{Ro}). 
Note that the center $\,\Gamma =\{\pm Id_2\}\,$  of $\,G\,$ acts   ineffectively on  $\,G^\C/K^\C\,$  
 and that $\,G/ \Gamma\,$ is isomorphic to
$\, SO_0(2,1)$. Then 
 the $\,SO_0(2,1)$-action on $\,D_\beta\,$ extends to
 $\,E(D_\beta) \,$ and the map $\, p\,$ is
equivariant, i.e.   $\,p \colon E(D_\beta) \to G^\C/K^\C \,$ is a
Stein, $\,SO_0(2,1)$-equivariant Riemann domain. 

By Theorem 7.5 in \cite{GeIa},   the map $\,p\,$ is necessarily injective.
Hence $\,E(D_\beta)\,$  coincides with the 
smallest, Stein,  $\,G$-invariant domain in $\,G^\C/K^\C\,$
containing $\,D_\beta$. Since each $\,D_\beta\,$  contains the Levi-flat orbit
$\,G\cdot w_1$, from the classification of Stein,  $\,G$-invariant domains in $\,G^\C/K^\C\,$ given in  Theorem 6.1 of \cite{GeIa},
  it follows that  the envelope of holomorphy
of $\,D_\beta\,$ is the domain $\,W$.

\end{proof}

\medskip
\nbigskip
{\bf Remark.}$\ $
Similarly one can show that $\,W\,$ is the envelope of holomorphy
of every element in the family of  
 hyperbolic, $\,G$-invariant subdomains of $\,W\,$
denoted in \cite{Is}, p.~22, by
$$\mathfrak D ^{(1)}_{s,t} \quad {\rm for} \quad -1\le s < 1 < t \le \infty,$$

 \nsmallskip
 where $\,s=-1\,$ and $\,t=\infty\,$ do not hold simultaneously.
 In terms of  diagram   (\ref{DIAGRAM1})
 the elements of the family $\,\mathfrak D ^{(1)}_{s,t}\,$ correspond  to the domains 
 $$ G \cdot (\ell_1((0,\alpha)) \cup   w_1 \cup   \ell_2((0,\beta)))
 \quad {\rm for} \quad 0 < \alpha \le1 \ \ {\rm and} \ \  0<\beta \le \infty \,.$$

 \bigskip

\end{document}